\documentclass[11pt]{article}
\usepackage{amssymb}
\usepackage[dvips]{graphics}
\usepackage{theorem}

\newtheorem{theorem}{{\sc Theorem}}[section]

\newtheorem{proposition}[theorem]{{\sc Proposition}}
\newtheorem{definition}{{\sc Definition}}[section]
\newtheorem{lemma}{{\sc Lemma}}[section]

\newtheorem{remark}[theorem]{{\sc Remark}}
\newtheorem{qquestion}[theorem]{{\sc Question}}
\newtheorem{pprobleme}[theorem]{{\sc Problem}}
\newtheorem{notation}{{\sc Notation}}

\newenvironment{proof}{{\sc Proof}}{\nolinebreak $\Box $}

\newenvironment{resume}{\small \begin{center} {\bf R\'esum\'e}
\end{center} \hspace{1cm} \begin{minipage}[t]{13cm} \hspace{.5cm} } 
{\end{minipage} \normalsize }

\begin{document}

\title{Combinatorially rigid simple polytopes with $d+3$ facets}
\author{Fr\'ed\'eric Bosio}
\date{\today}
\maketitle

\begin{resume}
  Nous classifions ici les polytopes simples combinatoirement rigides qui ont 
trois facettes de plus que leur dimension.
\end{resume}

\begin{abstract}
  We classify here combinatorially rigid simple polytopes with three facets 
more than their dimension.
\end{abstract}

\section{Introduction}

  We investigate combinatorial rigidity of simple polytopes. Different notions 
of rigidity of polytopes related to their combinatorics, and to toric theory, 
have been introduced. The first one is perhaps cohomological rigidity 
in~\cite{CPS}. It has the drawback to corcern only polytopes supporting 
quasitoric manifolds, so other notions have been introduced like combinatorial 
rigidity in~\cite{CK} (the one we consider here) or $B$-rigidity coming 
from~\cite{Bu}.

  It is not always easy to determine whether a polytope is rigid or not. 
Polygons are rigid. In dimension $3$, polyhedra without $3$ nor $4$-belt are 
$B$-rigid \cite{FMW}, whereas polyhedra with a $3$ or a $4$-belt are often 
nonrigid (if a polytope has a $3$-belt, it is rigid only when it separates the 
polytope in one part which is regular and the other is uniform, if it has a 
$4$-belt, it is nonrigid if, for instance, both parts have trivial isomorphism 
groups). In dimension $4$, a necessary condition for a polytope to be 
combinatorially rigid is that it is determined by the $1$-skeleton of its dual, 
i.e. by the "intersection graph" of its facets, which is already a strong 
property. Indeed, it still seems unknown if two polychora having the same 
intersection graph have diffeomorphic moment-angle manifolds. In higher 
dimension, even less is known.

  In this paper, we do not fix the dimension of the polytopes we study, but we 
require the difference between their number of facets and their dimension to 
equal $3$. We use the special combinatorial structure of these polytopes 
(theorem~\ref{structure}) to reduce the problem of their rigidity to a 
numerical problem, what eases its resolution.

\section{Statement of the problem}

\subsection{Recalls}

  A convex polytope is the convex hull of a (nonempty) finite set in a 
euclidean space. Its dimension is the one of the affine subpace it spans. Its 
faces are its intersections with supporting hyperplanes, its facets are its 
maximal faces (they are $1$-codimensional) and a $d$-dimensional polytope is 
called simple if all its vertices are contained in exactly $d$ facets (the 
minimum possible).

  The set of faces of a polytope is naturally partially ordered by inclusion 
and two polytopes are called isomorphic if their posets are. A combinatorial 
polytope is the poset of a polytope (indeed, it is determined by the inclusions 
of vertices into facets). Here, we only consider combinatorial polytopes, i.e. 
two polytopes with the same poset are identified.

  The following definition has been introduced in~\cite{CK}:

\begin{definition}
  A simple polytope $P$ is called {\em{combinatorially rigid}} if there is no 
other simple polytope with the same bigraded Betti numbers.
\end{definition}

  We investigate here the rigidity of $d$-dimensional simple polytopes with 
$d+3$ facets. The structure of these polytopes is known for a long 
(see~\cite{Gr}). I usually use the following version:

\begin{theorem}
\label{structure}
  A simple  $d$-dimensional polytope with $d+3$ facets is obtained by the 
operation of multiwedge over either the cube, or the dual cyclic polytope 
$C^* _{2k , 2k+3}$ for some $k \geq 1$.
\end{theorem}

  Thanks to this particular structure, we easily can get the Betti numbers of 
such a polytope.

  Let's before recall some basic facts about multiwedges~:

\begin{definition}
  Let $P$ be a simple polytope, its facets being numbered $F_1 ,...,F_n $. Let 
$M = (m_1 ,..., m_n )$ a $n$-tuple of natural numbers. Then the $M$-multiwedge 
over $P$ is the polytope with facets \newline
$F_{1,0} ,..., F_{1, m_1 } , F_{2,0} ,..., F_{2, m_2 } ,..., F_{n,0} ,..., 
F_{n, m_n }$ and so that a subset ${\cal F }$ of these facets determines a 
vertex of the multiwedge if:

- $\forall 1 \leq i \leq n$ there is at most one $j \leq m_n $ such that 
$F_{i,j} \notin {\cal F }$.

- The facets $F_i $ for whose each $F_{i,j}$ is in ${\cal F }$ determines a 
vertex of $P$.
\end{definition}

  Assume $P$ is a dual cyclic polytope $C^* _{2k , 2k+3}$ (for some basic facts 
about the structure of such a polytope, we refer to~\cite{Ga}). Its 
automorphism group is the dihedral group $D_{2(2k+3)}$, which is also the 
automorphism group of a $2k+3$-gon. We have a natural cyclic order (or dihedral 
order as it is unoriented) on the facets of $P$.

\begin{notation}
  We note $D$ the dihedral group $D_{4k+6}$.
\end{notation}

  The group $D$ acts naturally on the set of $(2k+3)$-tuples of integers (seen 
as maps from the sides of a $(2k+3)$-gon to ${\mathbb Z }$, and by restriction, 
on the set of $(2k+3)$-tuples of natural numbers.

\begin{remark}
Two $(2k+3)$-tuples give rise to the same combinatorial polytope if and only if 
they are in the same $D$-orbit.

  Indeed any simple polytope can be seen in a unique way as a multiwedge over a 
polytope which is not a wedge. In this constrcution, two facets are identified 
if any vertex in on either of them, and some vertex is on both. Hence two 
multiwedges on non-wedges $P_1 $ and $P_2 $ are combinatorially the same only 
if there is an isomorphism from $P_1 $ to $P_2 $ under which the $n$-tuples 
agree.
\end{remark}

  Let's now compute Betti numbers of our multiwedges. We use the following 
notation:

\begin{notation}
  Let $1 \leq i \leq 4k+6$ an integer. We note 
$\tilde{i} = \left\{\begin{array}{cl} 
i & \mbox{if $i \leq 2k+3$} \\
i - (2k+3) & \mbox {if $i > 2k+3$}
\end{array} \right. $
\end{notation}

  Obviously, in any case, $i$ and $\tilde{i}$ are congruent modulo $(2k+3)$.

  Recall that from Baskakov's formula \cite{Ba}, given two nonnegative integers 
$p$ and $q$, the Betti number $b^{p-q+1,2q}$ is the sum of the dimensions of 
the $p$-dimensional reduced homology groups of the sets obtained by union of 
$q$ facets of the polytope.

  The dual cyclic polytopes $C^* _{2k , 2k+3}$ have only two nontrivial nonzero 
Betti numbers, $b^{-1, 2k+2}$ and $b^{-2, 2k+4}$, whose value is $2k+3$. 
Indeed, apart from the emptyset and the wole polytope, the noncontratible 
unions of subsets of facets are the union of facets $F_1 , F_3 ,..., F_{2k+1}$, 
and their complements, up to cyclicity, and each of these unions has the 
homotopy type of a $k-1$-sphere.

  So we easily get the Betti numbers of the multiwedges we consider: For 
$1 \leq i \leq 2k+3$, the following union of facets:
$$\bigcup_{\begin{array}{c} 
0 \leq j \leq k \\ 0 \leq l \leq m_{\widetilde{i+2j}} 
\end{array}}
F_{\widetilde{i+2j} , l}$$
has the homotopy type of a sphere of dimension 
$(k-1) + \sum_{0 \leq j \leq k} m_{\widetilde{i+2j}}$, its complement has the 
homotopy type of a sphere of dimension 
$(k-1) + \left(\sum_{0 \leq j \leq k} m_{\widetilde{i+2j+1}}\right) + 
m_{\widetilde{i+2k+2}}$, whereas all other unions of facets, apart from the 
emptyset and the wole polytope, are contractible.

Each such subset contributes to $b^{-1, 2k + 2s_i + 2}$ where 
$s_i = \sum_{0 \leq j \leq k} m_{\widetilde{i+2j}}$, whereas, by Alexander 
duality for instance, its complement contributes to $b^{-2, 2k + 2s'_i + 4}$, 
where $s'_i $ is the sum of the other components of $M$.

  Remark that the equality of Betti numbers of two such multiwedges is 
equivalent to the equality of the former ones, which is equivalent to the fact 
that the sums $s_i $ thereabove take the same values the same number of times, 
independantly of their association with their indices.

\begin{definition}
  We say that two $(2k+3)$-tuples are {\em{$D$-equal}} if they belong to the 
same $D$-orbit.

  We say that two $(2k+3)$-tuples are {\em{equivalent}} if they give rise to 
multiwedges with the same Betti numbers.

  A $(2k+3)$-tuple is said rigid if the associated multiwedge is.
\end{definition}

\subsection{Sum-lists and numeric problem}

\subsubsection{Reformulation of the problem}

  We reformulate here the problem in more easier terms.

\begin{definition}
  Let $M = (m_1 ,..., m_{2k+3})$ a $(2k+3)$-tuple. For $1 \leq i \leq 2k+3$, we 
note like thereabove:
$$s_i = \sum_{j=0}^{k} m_{\widetilde{i+2j}}$$

  The list $L(M) = [s_1 ,..., s_{2k+3}]$ will be called the sum-list of this 
$(2k+3)$-tuple. We will note it $L$ if there is no ambiguity.
\end{definition}

  Let's remark that the sum-lists associated to two $D$-equal $(2k+3)$-tuples 
are $D$-equal.

  The foregoing tells us:

\begin{proposition}
  Two $(2k+3)$-tuples are $D$-equal if and only if their sum-lists are equal, 
up to dihedral order.

  Two $(2k+3)$-tuples are equivalent if and only if their sum-lists are equal, 
up to order.
\end{proposition}

  This reduces the problem of rigidity of polytopes with $d+3$ facets to a 
numerical problem: A $(2k+3)$-tuple is rigid if and only if any other 
$(2k+3)$-tuple giving rise to the same sum-list is $D$-equal to it.

  Indeed, the basic theory or circulant systems tells us that, given a 
permutation $L^{\sigma }$ of the sum-list $L$ of a $(2k+3)$-tuple $M$, we can 
find a $(2k+3)$-tuple $M^{\sigma }$ of possibly negative integers whose 
associated sum-list is $L^{\sigma }$.

  {\bf{Warning}}: The list $M^{\sigma }$ is not a permutation of the list $M$.

  More precisely, noting $S$ the sum of all the components of $M$, the 
components of $M^{\sigma }$ are the differences between $S$ and the sums of two 
consecutive (adjacent) elements of $L^{\sigma }$. Hence the fact that all 
components of $M^{\sigma }$ are nonnegative is equivalent to the fact that the 
greatest sum of two consecutive elements of $L^{\sigma }$ does not exceed $S$.

\begin{definition}
  A list $L^{\sigma }$ obtained from a permutation of the elements of the 
sum-list of a $(2k+3)$-tuple will be called a {\em{configuration}}.

  A configuration will be called {\em{admissible}} if the maximal sum of two 
consecutive elements of this configuration is not grater than $S$.
\end{definition}

  We now can reformulate the problem more simply:

\underline{Reformulation:} The $(2k+3)$-tuple $M$ is rigid if and only if any 
admissible configuration of $L$ is $D$-equal to $L$.

  This reformulation turns out to be helpful for solving the problem.

\subsubsection{Minimising a maximal sum}

  The problem we deal with concerns the minimisation of the maximal sum of 
two consecutive elements of a (cyclic) list. So we can ask: Given a list of 
numbers, what kinds of configurations minimise this maximal sum~?

\begin{definition}
  Let $L = [r_1 ,..., r_{2k+3}]$ a list of real numbers. We note:
$$K(L) = \max (r_{2k+3} + r_1 , \max_{1 \leq i \leq 2k+2} r_i + r_{i+1})$$
\end{definition}

  Clearly, $K(L) = K(L')$ if $L$ and $L'$ are $D$-equal.

\begin{notation}
  We consider a list $L$ of $2k+3$ real numbers.

  We order (by the usual real order) the elements of $L$, noting them in the 
following way: \newline
$x_0 \leq x_1 \leq ... \leq x_{k+1} \leq y_{k+1} \leq ... \leq y_1 $.

  We also order the different values of the elements the sum list~: Assume the 
sum-list contains $r+1$ different values. We order them increasingly 
$u_0 < u_1 <...< u_r $ or decreasingly $v_0 > v_1 >...> v_r $ (so 
$u_i = v_{r-i}$). We can also notice $u_0 = x_0 $ and $v_0 = y_1 $.

  If $v$ is a value of the list, we note $m_v $ its multiplicity in the 
sum-list (i.e. the number of times it appears in it). The multiplicity 
$m_{v_0 }$ of the greatest value of the sum-list will be noted $m$.

  We also note $K = \max_{1 \leq i \leq k+1} (x_i + y_i )$.
\end{notation}

\begin{proposition}
  Let $L$ a list of $2k+3$ real numbers.

  Now, consider a configuration $L^{\sigma }$ of the elements of $L$.

  Then, we have $K(L) \geq K$ and there is a configuration $L^{\sigma }$ so 
that $K(L^{\sigma }) = K$.
\end{proposition}

  We even will provide many ways to configurate the elements so that 
$K(L^{\sigma }) = K$.

\begin{proof}
  Let's first prove $K(L) \geq K$. Let $1 \leq i \leq k+1$ and consider the 
elements $y_1 ,..., y_i $. Then either two of them are adjacent, in which case 
$K' \geq 2y_i \geq x_i + y_i $, or they globally have at least $i+1$ 
neighbours, so one of these neighbours is at least $x_i $, and the sum of the 
two forementioned neighbours is at least $x_i + y_i $. So we have 
$K' \geq x_i + y_i $ in any case, which implies $K(L) \geq K$.

  Let's now describe configurations for whose $K(L^{\sigma }) = K$. The 
proof thereabove suggests to find configurations for which, given 
$1 \leq i \leq k+1$, the only neighbours of $y_1 ,..., y_i $ are 
$x_0 , x_1 ,..., x_i $. Such configurations exist. On the first step, we place 
$y_1 $ surrounded by $x_0 $ (on its right) and $x_1 $ (on its left). For the 
moment, $x_0 $ and $x_1 $ have only one neighbour. They are said available. On 
the second step, we put $y_2 $ adjacent to an available element ($x_0 $ or 
$x_1 $) and $x_2 $ adjacent to $y_2 $. Now the available elements are $x_2 $ 
and the element $x_0 $ or $x_1 $ which is not adjacent to $y_2 $. We can 
continue this process. At step $i$, we put $y_i $ adjacent to an available 
element and put $x_i $ adjacent to $y_i $. We stop when all elements have been 
put in the configuration, i.e. after step $k+1$.

  For such a configuration, we see that two consecutive elements can be:

\begin{itemize}
\item
$y_1 + x_0 \leq y_1 + x_1 \leq K$.

\item
$y_i + x_i \leq K$.

\item 
The sum of $y_i $ with an element which was available at step $i$, so this sum 
is $\leq y_i + x_{i-1} \leq K$

\item
The sum of $x_{k+1}$ with an element which was available at the last step, so 
this sum is $\leq x_{k+1} + x_k \leq K$.
\end{itemize}

  In any case, the sum if at most $K$, so $K(L^{\sigma }) = K$ in any 
configuration we have constructed thereabove.
\end{proof}

\begin{remark}
  We notice that if we consider a quite generic sum-list, many non 
$D$-equal configurations will produce a non greater $K(L^{\sigma })$. So a 
generic polytope with $d+3$ facets won't be rigid.
\end{remark}

  Consider the configurations we have constructed thereabove. At step $i$, we 
say we have made the choice $l$ (left) if, for the cyclic order we get at the 
end, $y_1 , x_1 , y_i , x_i $ appear in this order, and we say we have made the 
choice $r$ (right) if, for the cyclic order we get at the end, 
$y_1 , x_0 , y_i , x_i $ appear in this order. Such a configuration will be 
encoded by the list of choices.

  Finally, we call {\em{standard}} the configuration in which we make only 
choices $l$, so the following one:

\includegraphics{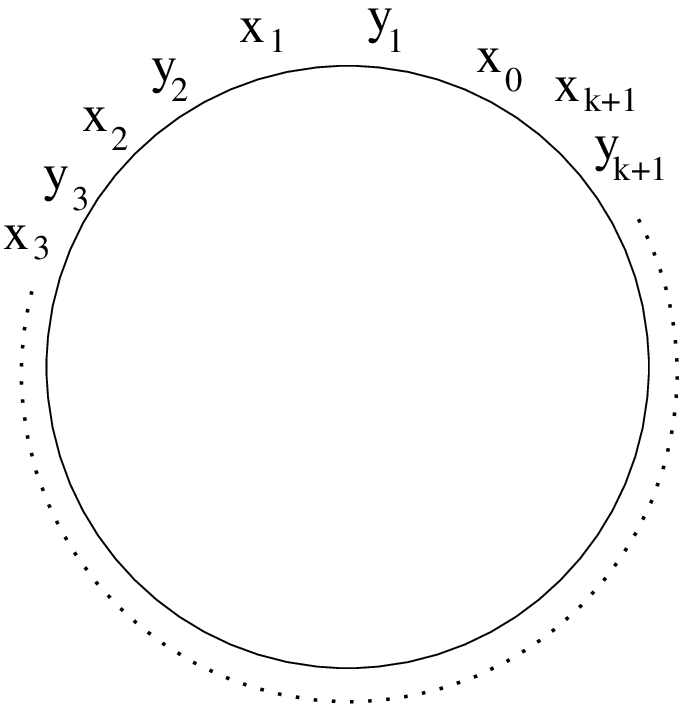}

  We will use the following terminology:

\begin{definition}
  Consider a list of $(2k+3)$ numbers. Two elements are called {\em{adjacent}} 
or {\em{consecutive}} if their indices differ from $1$ or $2k+2$. If $l$ is an 
element of the list, an element $l'$ adjacent to $l$ is also called a 
{\em{neighbour}} of $l$. 

  Consider a list of $(2k+3)$ numbers. Then, the elements of a sublist of 
$2 \leq j \leq k+1$ terms are called {\em{jump-adjacent}} if
there is $i$ such that the indices of these elements are 
$i , \tilde{i+2} , \tilde{i+4} ,..., \tilde{i + 2j - 2}$.
\end{definition}

  We immediately check that these notions are preserved by $D$. We remark also 
that the $j \geq 2$ elements of a sublist of a list are jump-adjacent if and 
only if no two of them are adjacent but they globally have only $j+1$ 
neighbours.

\section{Determination of rigid polytopes}

  We now classify rigid  polytopes with $d+3$-facets. We can already state our 
main theorem:

\begin{theorem}
  If a polytope is a product of three simplices, then it is rigid.

  The rigid multiwedges over the pentagon are given by the following $5$-tuples:

\begin{enumerate}
\item
$(a,a,b,b,b)$.

\item
$(a,b,b,c,c)$ with $a>b+c$.

\item
$(a + \lambda , a , b , b + \lambda , c)$ where
$c < \lambda $.

\item
$(a,b,c,d,d)$ with $\min (a,c) > b+d$.
\end{enumerate}

  A multiwedge over the dual cyclic polytope $C^* _{2k , 2k+3}$, where 
$k \geq 2$ is rigid if and only if it is given by a $(2k+3)$-tuple of the 
following form:

\begin{enumerate}
\item
$(a,a,b,b,...,b)$.

\item
$(a,b,b,...,b)$ where $a>2b$.

\item
$(a,b,b,...,b,a,c,c,...,c)$ where $a > b+c$.

\item
$(a,b,c,d,...,d)$ where $\min (a,c) > b+d$.

\item
$(a,b,b,c,c,a,d,...,d)$ where $a > \max (b,c) + d$

\item
$(a + \lambda , a, b , b + \lambda , c , c , ... , c)$ where $\lambda > c$.

\item
$(a + \mu , a , a + \lambda , b , ... , b , a + \lambda + \mu , c ,..., c)$ 
where $\lambda > b$ and $\mu > c$.

\end{enumerate}
\end{theorem}

  As we mentioned earlier, the reformulation allows an easier way of solving 
the problem. The theorem will  be a direct consequence of the following result:

\begin{proposition}

  We give here the list of rigid sum-lists (i.e. the sum-lists of rigid 
$(2k+3)$-tuples). Given a list $[\lambda _1 ,..., \lambda _{2k+3}]$ of $2k+3$ 
elements, call $S = \frac{1}{k+1} \sum_{i=1}^{2k+3} \lambda _i $. Then, a list 
of $2k+3$ elements is rigid if and only if all its elements are nonnegative, 
the sum of any two adjacent elements does not exceed $S$, and it is $D$-equal 
to a member of one of the following families:

\begin{enumerate}
\item All elements but at most one have the same value.

\item $[u_0 , v_0 , u_1 , \underbrace{u ,..., u}_{2k \mbox{ times }}]$, where 
$v_0 + u > S$.

\item $[\underbrace{u_0 , v_0 ,..., u_0 , v_0 }_{(k+1) \mbox{ times }}, u_0 ]$, 
where $2v_0 > S$.

\item 
$[\underbrace{u_0 , v_0 ,..., u_0 , v_0 }_{k \mbox{ times }}, u_0 , u , v]$, 
where $v_0 + \min(u,v) > S$.

\item For $k>1$
$[u_0 , v_0 , u_0 , v_1 , u_0 , \underbrace{u ,..., u}_{2k-2 \mbox{ times }}]$, 
where $v_1 + u > S$.

\item $[u_0 , v_0 , u_0 , v_1 , u_0 ]$, where $v_0 + v_1 > 3 u_0 $.

\item $[\underbrace{u_0 , v_0 ,..., u_0 , v_0 }_{m<k \mbox{ times }} , u_0 , 
\underbrace{u ,..., u}_{2(k-m+1) \mbox{ times }}]$, where $v_0 + u > S$.

\item $[\underbrace{u_0 , v_0 ,..., u_0 , v_0 }_{m<k \mbox{ times }} , u_0 , 
v_1 , u_1 , \underbrace{u ,..., u}_{2(k-m) \mbox{ times }}]$, where 
$v_0 + u_1 > S$ and $v_1 + u > S$.

\end{enumerate}
\end{proposition}

  Before proving the proposition, we establish the correspondance between the 
rigid $(2k+3)$-tuples of the theorem and the sum-lists of the proposition. Each 
case is a straightforward computation.

  The case 1 of the proposition corresponds to both cases 1 of the theorem.

  The case 2 of the proposition corresponds, if $k=1$, to case 3 of the 
theorem, and if $k>1$, to case 6 of the theorem.

  The case 3 of the proposition corresponds, if $k=1$, to case 2 of the theorem 
where $b=c$ or if $k>1$, to case 2 of the theorem.

  The case 4 of the proposition corresponds to both cases 4 of the theorem.

  The case 5 of the proposition corresponds to $k>1$, case 5 of the theorem.

  The case 6 of the proposition corresponds to $k=1$, case 2 of the theorem 
where $b \neq c$.

  The case 7 of the proposition corresponds to $k>1$, case 3 of the theorem.

  The case 8 of the proposition corresponds to $k>1$, case 7 of the theorem.

  We then prove the proposition~:

\begin{proof}
  Let's first check that all cases or the proposition are actually rigid.

  In case 1, it is clear.

  In case 2, the element $u_0 $ cannot have $v$ as neighbour, so his neighbours 
must be $u_0 $ and $u_1 $. Then, up to $D$, these three elements can be put as 
in the standard case. As all others are equal, we are done.

  In case 3, there are $k+1$ occurences of $v_0 $ and no two of them can be 
consecutive. The only possibility is that they are jump-adjacent, so a unique 
configuration up to $D$ (and even up to the cyclic group).

  In case 4, the occurences of $v_0 $ cannot have any other neighbour than 
$u_0 $. As $u_0 $ appears only one more time than $v_0 $, the occurences of 
$v_0 $ must be jump-adjacent, and all their neighbours be $u_0 $. So they are 
up to $D$ in standard position. Only two elements are remaining, but exchanging 
them yields a symmetry of the configuration. So any two admissible 
configurations are $D$-equal.

  In case 5, both $v_0 $ and $v_1 $ must have their neighbours equalling 
$u_0 $. As there are only three occurences of $u_0 $, they must be 
jump-adjacent. So the positions of the occurences of $u_0 $, $v_0 $ and $v_1 $ 
are fixed up to $D$-symmetry. As all other coefficients are equal, all 
admissible configurations are $D$-equal.

  In case 6, as three elements have the same value, there are only two 
$D$-orbits of configurations with the elements of the list, the one in which 
$v_0 $ and $v_1 $ are adjacent, the one in which they are not. Due to the 
inequality, only the latter is admissible, so the sum-list is actually rigid.

  In case 7, the occurences of $v_0 $ cannot have any other neighbour 
than $u_0 $. As $u_0 $ appears only one more time than $v_0 $, the occurences 
of $v_0 $ must be jump-adjacent, and all their neighbours be $u_0 $. So they 
are up to $D$ in standard position. As all other coefficients are equal, all 
admissible configurations are $D$-equal.

  In case 8, the occurences of $v_0 $ cannot have any other neighbour 
than $u_0 $. As $u_0 $ appears only one more time than $v_0 $, the occurences 
of $v_0 $ must be jump-adjacent, and all their neighbours be $u_0 $. So they 
are up to $D$ in standard position. The only occurence of $v_1 $ cannot have 
any other neighbour than $u_0 $ and $u_1 $. As it can have at most one 
neighbour of each value, it must have a neighbour equal to $u_0 $ and the other 
equal to $u_1 $. So the positions of $u_1 $ and $v_1 $ are also fixed up to 
$D$. As all other coefficients are equal, all admissible configurations are 
$D$-equal.

  This already proves the rigidity of the announced sum-lists.

  We now have to prove the converse, i.e. that our list contains all the rigid 
sum-lists.

  We now fix a $(2k+3)$-tuple $M$ and define all as above ($L$, $S$, etc...). 
We dinstinguish several cases:

\vskip 3mm

\underline{First case: Completely rigid case.}

\begin{definition}
  A $(2k+3)$-tuple is called {\em{completely rigid}} if all configurations are 
$D$-equal.
\end{definition}

  Clearly, a completely rigid $(2k+3)$-tuple is rigid.

  As the dihedral group $D_6 $ is the same as the permutation group 
$\mathfrak{S}_3 $, any triple of nonnegative integers is completely rigid. So 
il the sequel, we always assume $k \geq 1$.

 Let's consider a partition of $L$ into two sublists having no common value, 
namely $L_1 $ and $L_2$. There is a configuration for which the elements of 
$L_1 $ are consecutive. So these elements must be consecutive in any 
configuration. If there are at least two elements in both $L_1 $ and $L_2$, 
then there is a configuration where an element of $L_1 $ has its two neighbours 
in $L_2 $ and the elements of $L_1 $ are then not consecutive. We then cannot 
be in the completely rigid case.

  So every value is taken either only once or at least at least $2k+2$ times. 
And, as $k \geq 1$, we cannot have two values taken only once. So all elements 
of the sum-list but at most one have the same value.

  Conversely, if all elements of the sum-list but at most one are equal, the 
configuration is clearly completely rigid.

\vskip 3mm

\underline{Second case: Non unique small value.}

\begin{definition}
  A value in the sum-list is called {\em{small}} if its sum with $v_0 $ does 
not exceed $S$.
\end{definition}

\begin{remark}
  Any permutation of occurences of small values from an admissible 
configuration yields another admissible configuration.

  Clearly, the sums of two consecutive values in the new configuration are 
either a sum of a small element and an element of the sum-list, or are sum of 
two consecutive values in the original admissible configuration. In both cases, 
they do not exceed $S$.
\end{remark}

\begin{lemma}
  Assume there are at least two small values and the configuration is not 
completely rigid. Then there are exactly two small values, each has 
multiplicity $1$ and $v_0 $ has multiplicity $1$ too.
\end{lemma}

\begin{proof}
  Thanks to the precedent remark, the two small extremal elements of the list 
may be any small value. So any sublist of two extremal values must be the same 
up to permutation. So, if there are two different values in the sum-list, then 
any sublist with two elements of the list of small values must contain both 
these two elements. So the list of small values must come down to these two 
elements. Moreover, as the neighbours of all occurences of $v_0 $ must be 
small, the multiplicity of $v_0 $ is at most $2-1$, so must be $1$.
\end{proof}

  We now settle this case. Assume there are (at least) two small values and the 
configuration is not completely rigid. Then $u_2 = v_1 $, i.e. all values 
except $u_0 $, $u_1 $ and $v_0 $ are equal, and so we are in case 2 of the 
proposition.

  Indeed, in this case, if we fix $v_0 $ and $u_0 $ as in the standard case, 
the only admissible configuration is the standard one. In particular, the 
$r^k $ configuration is identical to the standard one, so 
$y_2 = x_{k+1} \leq y_{k+1} = x_2 $. This proves the claim. We just need the 
inequality $v_0 + u_2 > S$ to guarantee that $u_0 $ and $u_1 $ are the only 
small value.

\vskip 3mm

\underline{Third case: $m=k+1$.}

  We here assume $m=k+1$. This case is easy. In this case, the configuration 
cannot be completely rigid and $m>1$ so we have a unique small value $u_0 $ so 
$m_{u_0 } \geq m+1 = k+2$. Also, $m_{u_0 } \leq 2k+3 - m = k+2$. Then we only 
have two values, $u_0 $ appearing $k+2$ times and $v_0 $ appearing $k+1$ times. 
This corresponds to case 3 of the proposition. The only thing we need here to 
guarantee rigidity is that $v_0 $ is not a small value, i.e. the inequality 
$2 v_0 > S$.

\vskip 3mm

\underline{Fourth case: $m = k > 1$.}

  We here assume $m = k > 1$. This case is not much more difficult. In this 
case, the configuration cannot be completely rigid and $m>1$ so we have a 
unique small value $u_0 $ so $m_{u_0 } \geq m+1 = k+1$. Hence, apart from 
$u_0 $ and $v_0 $, there are am most two others elements in the sum-list. 
Indeed, we must have $m_{u_0 } = k+1$. Had we $m_{u_0 } > k+1$ would we have 
$x_{k+1} = u_0 $. From the standard configuration, putting it between $y_0 $ 
and $x_1 $ would yield another admissible configuration which is not $D$-equal 
to the standard one as the occurences of $v_0 $ would not be jump-adjacent.

  So we are in the case 4 of the proposition. The only inequality we then need 
to guarantee rigidity is $v_0 + u_1 > S$.

\vskip 3mm

\underline{Fifth case: $m_{u_0 } > m+1$.}

  We assume here the configuration is not completely rigid and 
$m_{u_0 } > m+1$. We see immediately we're not in the second case so $u_0 $ is 
the only small value.

\begin{lemma}
  Under these hypotheses, we must have $m_{u_0 } = 3$ and 
$m = m_{v_1 } = 1$.
\end{lemma}

\begin{proof}
  There must be some occurence of $u_0 $ which is not adjacent to any of 
$v_0 $. If we consider the standard configuration, any permutation of the 
elements of the list both of whose neighbours are $u_0 $ is admissible. So 
permuting $y_m = v_0 $ with $y_{m+1} = v_1 $ must give a $D$-equal 
configuration. We then cannot have $m>1$, as in this case, the occurences of 
$v_0 $ would no more be jump-adjacent.

  We neither can have $m_{v_1 } > 1$. Indeed, the neighbours of the neighbours 
of $v_0 $ are fixed by $D$, so, as $m=1$ these two neighbours are only 
$y_2 = v_1 $ and $x_{k+1}$ in the standard configuration. They are $y_2 = v_1 $ 
and $y_3 $ after the transposition of $y_1 $ and $y_2 $. We then must have 
$y_3 = x_{k+1}$. If we had $m_{v_1 } > 1$, we would have 
$y_3 = y_2 = x_{k+1} = y_{k+1}$ and so $2 v_1 = x_{k+1} + y_{k+1} \leq S$. So 
any permutation in which $v_0 $ has only $u_0 $ as neighbours would be 
admissible, so $D$-equal to the standard one. But in some such configuration 
all the occurences of $v_1 $ are consecutive, and not in the standard one (we 
cannot have $v_1 = u_0 $ as the configuration is assumed not completely 
rigid). So the hypothese $m_{v_1 } > 1$ leads to a contradiction.
\end{proof}

  We then settle this case. If $k=1 $, we are in case 6 of the proposition. The 
only thing we need to guarantee rigidity is that $v_0 + v_1 $ must be gretaer 
than $S$, which is equivalent to $v_0 + v_1 > 3 u_0 $.

  Else, all other elements apart from the five equalling $u_0 $, $v_0 $ or 
$v_1 $ are equal. Indeed, any admissible configuration fixing $u_0 $ and $u_1 $ 
must be the standard one. This is the case of the $lr^{k-1}$ configuration. So 
we have $y_3 = x_{k+1} \leq y_{k+1} = x_3 $. This proves the claim. This 
corresponds to case 5 of the proposition. We only need the inequality 
$v_1 + u_1 > S$ to avoid adverse admissible configurations.

\vskip 3mm

\underline{Sixth case: Remaining cases.}

  We here have $m \leq k-1$ occurences of $v_0 $ and $m+1$ occurences of 
$u_0 $, the unique small value. There are then exactly two occurences of $u_0 $ 
that have a neighbour other than $v_0 $. In a given configuration, only two 
elements of $D$ act by globally preserving these two occurences of $u_0 $, the 
identity and a symmetry. Hence any admissible configuration putting the $v_0 $ 
at the same positions as the standard one must be either the standard one 
itself or its symmetric exchanging these two occurences of $u_0 $.

 In the standard configuration, these neighbours are $y_{m+1} = v_1 $ and 
$x_{k+1}$.

  If they are equal, they also are equal to $y_{k+1}$, so 
$2 v_1 = x_{k+1} + y_{k+1} \leq S$ and all configuration obtained from the 
standard one by permuting the elements other than $u_0 $ and $v_0 $ must be 
$D$-equal. As there are at least $2k+3 - (2k-1) = 4$ such elements, they all 
must be equal. We then are in case 7 of the proposition. The only thing we need 
to guarantee rigidity is $v_0 + u > S$ where $u$ is the only value but 
$u_0 $ and $v_0 $.

  Assume now $v_1 \neq x_{k+1}$. Then the configuration $l^{m-1}rl^{k-m}$ 
cannot be equal to the standard one, so must be its forementioned symmetric. 
We then get $y_{m+2} = x_{k+1} \leq y_{k+1} = x_{m+2}$. So we are in case 8 of 
the proposition.

  The exchange of $y_{m+1}$ with $x_{m+1}$ cannot be admissible, which gives 
$v_1 + u_2 > S$. Neither can be the exchange of $x_{m+1}$ with $x_m $, which 
gives $v_0 + u_1 > S$.

  We have considered all the possibilities. So the proposition is proved.
\end{proof}

\section{Final remarks}

  We give here some remarks and comments about our result.

\paragraph{Other rigidities}
  We have considered here the notion of combinatorial rigidity, as it was the 
easiest to compute. The moment-angle manifold associated to such a polytope is 
the connected sum of sphere products (this is proved in \cite{LdMV} although 
authors were not aware they dealt with such moment-angle manifolds), so is 
determined by it Betti numbers. The $B$-rigidity is not in this case really 
different from the combinatorial rigidity. The cohomological rigidity assumes 
the existence of quasitoric manifolds over the considered polytopes. Indeed, 
for $k$ large, there is no quasitoric manifold over the dual cyclic polytope 
$C^* _{2k , 2k+3}$ \cite{Ha}, hence the notion of cohomological rigidity seems 
less meaningful.

\paragraph{Smallest nonrigid polytope}
\begin{remark}
  We can see that $(m_1 ,..., m_n )$ is rigid as long as $\sum _i m_i \leq 2$. 
It is obvious for $0$ and $1$, and if $\sum _i m_i = 2$, the number of $s_i $ 
equalling $2$ (or $0$) indicates the relative positions (up to $D$) of the 
facets on which the wedges are performed.

  So the smallest nonrigid polytopes appear when $k=1$ and $\sum _i m_i = 3$, 
i.e. for $5$-dimensional polytopes with eight facets. Up to dimension $4$, all 
simple polytopes with $d+3$ facets are rigid.

  The $5$-tuples $(2,1,0,0,0)$ and $(1,1,0,1,0)$ are equivalent (both give the 
sum-list $(0,1,1,2,2)$ up to order), but they are not $D$-equal.

  Both corresponding polytopes are obtained by a single wedge over the polytope 
given by the $5$-tuple $(1,1,0,0,0)$. Notice that for the first-mentioned 
polytope, the facet on which this last wedge is made is a pentagonal book, 
whereas for the second-mentioned polytope, it is a cube. This provides an 
example of two wedges over nonisomorphic facets of the same polytope that give 
diffeomorphic moment-angle manifolds.
\end{remark}

\paragraph{Rigidity locus}
  It could also be interesting to look at rigid $(2k+3)$-tuples is some parts 
of the space ${\mathbb N }^{2k+3}$. We just look here around the diagonal, i.e. 
when the components $m_i $ are close to each other. A direct corollary of our 
theorem is:

\begin{remark}
  Consider a rigid $(2k+3)$-tuple $(m_1 ,..., m_{2k+3})$ so that 
$\max_i m_i \leq 2 \min_i m_i $.

  Then $(m_1 ,..., m_{2k+3})$ is completely rigid.
\end{remark}

\paragraph{Semi-projectivization}
  Another simple remark is that the rigidity is preserved by multiplication, 
i.e. given a $(2k+3)$-tuple $(m_1 ,..., m_{2k+3})$ and an integer $n \geq 1$, 
then $(m_1 ,..., m_{2k+3})$ is rigid if and only if 
$(n \cdot m_1 ,..., n \cdot m_{2k+3})$ is.

  So we can, in some sense, look at a "semi-projectivization" of the problem. 
We consider the set of $(2k+3)$ nonnegative numbers whose sum equals $1$, i.e. 
the standard simplex (of dimension $(2k+2)$). The permutation group 
${\mathfrak S }_{2k+3}$ naturally acts on this simplex. The dihedral group 
$D$ is a subgroup of ${\mathfrak S }_{2k+3}$ so also naturally acts on this 
simplex. Moreover, its action preserves the following subset:
$$Z_{2k+3} = \{ (x_1 ,..., x_{2k+3}) , x_i \geq 0 , \sum_i x_i = 1, 
x_i + x_{\widetilde{i+1}} \leq \frac{1}{k} \} $$

  Then, up to some kind of semi-projectivization, the rigid sum-lists 
correspond to the $(2k+3)$-tuples $(x_1 ,..., x_{2k+3})$ of $Z_{2k+3}$ for 
whose the intersection of their ${\mathfrak S }_{2k+3}$-orbit with $Z_{2k+3}$ 
is reduced to their $D$-orbit.

  This leads to other problems of the same type~: We can consider the 
odd-dimensional analogue, change the value $\frac{1}{k}$ in the definition of 
$Z_{2k+3}$, consider the sum of more than two consecutive terms of the list ... 
all these problems may deserve consideration.

\paragraph{One more facet}
  Finally, let's say a few word on (very particular considerations) about the 
rigidity problem for simple polytopes with $d+4$ facets. In an unpublished 
work, the author has classified the pairs of multiwedges over the hexagon 
that have the same Betti numbers. Indeed, there are basically four families of 
pairs of $6$-tuples, each family having five parameters, that give multiwedges 
over the hexagon with the same Betti numbers. So, generically, given a 
multiwedge over the hexagon, we cannot find another one with the same Betti 
numbers. But this does not allow us to conclude uncautiously that such a 
polytope is rigid, as other polytopes than such multiwedges might have the same 
Betti numbers (indeed, there is no known structure theorem for simple polytopes 
with $d+4$ facets). Let's give $4$-dimensional examples:

  Consider the biwedges over the hexagon given by the $6$-tuples 
$(1,0,1,0,0,0)$ and $(1,0,0,1,0,0)$. They have the same Betti numbers, as well 
as three other polytopes. Indeed the three polytopes obtained by vertex 
truncation on the $(1,1,0,0,0)$-wedge over the pentagon (including the 
$(1,0,1,0,0,0)$-multiwedge over the hexagon), and the two polytopes obtained 
by connected sum of two copies of the product of two triangles (including the 
$(1,0,0,1,0,0)$-multiwedge over the hexagon) have the same Betti numbers.

  Consider the biwedge over the hexagon given by the $6$-tuple $(1,1,0,0,0,0)$. 
No other multiwedge over the hexagon has the same Betti numbers but three other 
polytopes have. Indeed, all these polytopes are obtained by truncating and edge 
of the $(1,1,0,0,0)$-wedge over the pentagon (notice there is another polytope, 
with different Betti numbers obtained by this construction, namely if we 
cut off a "horizontal" edge of the unique cubical facet of this 
$(1,1,0,0,0)$-wedge).

  So the rigidity of a generic multiwedge over the hexagon is unknown but the 
precedent examples may suggest its falsity.

{\footnotesize {Bosio Fr\'ed\'eric \\
UMR 7348 \\
UFR Sciences SP2MI \\
Teleport 2 \\
Boulevard Marie et Pierre Curie \\
BP 30179 \\
86962 Futuroscope Chasseneuil CEDEX 

e-mail~: bosio@math.univ-poitiers.fr}}


\begin{thebibliography}{99}
\bibitem[Ba]{Ba} I. Baskakov {\em Cohomology of K-powers of spaces and
the combinatorics of simplicial divisions} \newline
Russian math. surveys 57 (2002), no. 5, p. 989-990

\bibitem[B-P1]{BP1} V. Buchstaber and T. Panov {\em Torus actions and
their applications in topology and combinatorics} \newline
University lectures series, 24. Amer. math. soc., providence, RI, 2002.

\bibitem[B-P2]{BP2} V. Buchstaber and T. Panov {\em Toric topology} \newline
Mathematical Surveys and Monographs, vol. 204, Amer. math. soc., 2015.

\bibitem[Bu]{Bu} V. Buchstaber {\em Lectures on toric topology. Lecture
notes of 'Toric topology workshop: KAIST 2008'} \newline
Trends in Math. 10 (2008), no. 1, p. 1-64

\bibitem[CK]{CK} S. Choi and J.S. Kim {\em Combinatorial rigidity of
$3$-dimensional simplicial polytopes} \newline
Int. Math. Res. Not. IMRN. 2011, no.8, p. 1935-1951

\bibitem[CMS]{CMS} S. Choi, M. Masuda and D.Y. Suh {\em Rigidity
problems in toric topology: a survey} \newline
Tr. Mat. Inst. Steklova, 2011, Volume 275, p. 188–201

\bibitem[CPS]{CPS} S. Choi, T. Panov and Y. Suh {\em Toric
cohomological rigidity of simple convex polytopes} \newline
Journal of London Math. Soc. II, Ser. 82 (2010), p. 343-360

\bibitem[D-J]{DJ} M. Davis and T. Januszkiewick {\em Convex polytopes,
Coxeter orbifolds and torus actions} \newline
Duke Math. Journal, vol 62 (1991), no. 2, p. 417-451

\bibitem[FMW]{FMW} F. Fan, J. Ma and X. Wang {\em $B$-rigidity of flag
$2$-speres without $4$-belt} \newline
arXiv:1511.03624

\bibitem[Ga]{Ga} D. Gale. {\em Neighborly and cyclic polytopes} \newline
 Proc. Sympos. Pure Math., Vol. VII, p. 225-232. Amer.Math.Soc.,
 Providence, R.I., 1963.

\bibitem[Gr]{Gr} B. Gr\"unbaum {\em Convex polytopes}, 2nde \'edition \newline
Graduate texts in mathematics (221). Spinger, 2003.

\bibitem[Ha]{Ha} S. Hasui
{\em On the classification of quasitoric manifolds over dual cyclic 
polytopes} \newline
Algebraic \& Geometric Topology 15 (2015), p. 1387-1437

\bibitem[LdM-V]{LdMV} S. L\'opez de Medrano and A. Verjovsky {\em A
new family of complex, compact, non-symplectic manifolds} \newline
Bol. Soc. Brasil. Mat. (N.S.) 28 (1997), no. 2, p. 253-269

\end{thebibliography}
\end{document}